\documentclass[leqno,a4paper,12pt]{article}
\usepackage{amsmath,amsthm,bbm}
\usepackage[sans]{dsfont}

\newtheorem{Thm}{\indent Theorem}[section]
\newtheorem{Prop}[Thm]{\indent Proposition}
\newtheorem{Lem}[Thm]{\indent Lemma}
\newtheorem{Cor}[Thm]{\indent Corollary}

\theoremstyle{definition}
\newtheorem{Def}[Thm]{\indent Definition}
\newtheorem{Rem}[Thm]{\indent Remark}

\def\qed{{\hskip0pt\unskip\unskip\nobreak\hfil\penalty50
          \hskip1em\hbox{}\nobreak\hfil
          {\bf q.e.d.}%
          \parfillskip=0pt\finalhyphendemerits=0
          \par}\medskip}

\newenvironment{Proof}
               {{\it Proof.}\quad}
               {\qed}

\newcommand{\Prime}{\kern3\fontdimen1\font$'$\kern-7\fontdimen1\font}

\long\def\forget#1{}

\long\def\beginSIDEREMARK#1\endSIDEREMARK
    {{\par\bigskip\advance\leftskip by 2cm
                  \advance\rightskip by -2cm\noindent
      {\bf Our own side remark:} #1
      \par\bigskip\noindent}}

\long\def\beginFORGET#1\endFORGET{#1}
\long\def\beginFORGET#1\endFORGET{}

\def\?{\ ???\ \immediate\write16{}%
\immediate\write16{Warning: There was still a question mark . . . }%
\immediate\write16{}}

\usepackage{amsmath}

\usepackage{exscale}
\usepackage{amssymb}

\input cyracc.def
\font\tencyr=wncyr6
\def\cyr{\tencyr\cyracc}
\newcommand{\cyrb}{{\cyr B}}

\newcommand{\BB}{{\mathbb{B}}}
\newcommand{\BC}{{\mathbb{C}}}

\newcommand{\BQ}{{\mathbb{Q}}}

\newcommand{\BZ}{{\mathbb{Z}}}

\newcommand{\Fc}{{\mathfrak{c}}}

\newcommand{\Fn}{{\mathfrak{n}}}

\newcommand{\FA}{{\mathfrak{A}}}

\newcommand{\FS}{{\mathfrak{S}}}

\newcommand{\CC}{{\cal C}}
\newcommand{\CD}{{\cal D}}

\newcommand{\CH}{{\cal H}}

\newcommand{\CK}{{\cal K}}
\newcommand{\CL}{{\cal L}}

\newcommand{\Spec}{\mathop{{\bf Spec}}\nolimits}

\newcommand{\rad}{\mathop{{\rm rad}}\nolimits}

\newcommand{\GL}{\mathop{\rm GL}\nolimits}

\newcommand{\Gr}{\mathop{\rm Gr}\nolimits}

\newcommand{\Hom}{\mathop{\rm Hom}\nolimits}

\newcommand{\loccit}{[loc.$\;$cit.]}

\def\halb{\frac{1}{2}}

\newbox\mybox
\def\arrover#1{\mathrel{
       \setbox\mybox=\hbox spread 1.4em{\hfil$\scriptstyle#1$\hfil}
       \vbox{\offinterlineskip\copy\mybox
             \hbox to\wd\mybox{\rightarrowfill}}}}
\def\larrover#1{\mathrel{
       \setbox\mybox=\hbox spread 1.4em{\hfil$\scriptstyle#1$\hfil}
       \vbox{\offinterlineskip\copy\mybox
             \hbox to\wd\mybox{\leftarrowfill}}}}

\def\ontoover#1{\mathrel{
       \setbox\mybox=\hbox spread 1.4em{\hfil$\scriptstyle#1$\hfil}
       \vbox{\offinterlineskip\copy\mybox
             \hbox to\wd\mybox{\rightarrowfill\hskip-2.8mm
                               $\rightarrow$}}}}
\def\leftontoover#1{\mathrel{
       \setbox\mybox=\hbox spread 1.4em{\hfil$\scriptstyle#1$\hfil}
       \vbox{\offinterlineskip\copy\mybox
             \hbox to\wd\mybox{$\leftarrow$\hskip-2.8mm
                               \leftarrowfill}}}}
\def\longto{\longrightarrow}
\def\into{\hookrightarrow}

\def\longinto{\lhook\joinrel\longrightarrow}

\usepackage[curve,matrix,arrow,cmtip]{xy}
\NoComputerModernTips

\def\myxymessage{\def\messagetext
   {Here an xy-pic diagram was omitted to speed up compilation . . . }
   \immediate\write16{\messagetext}
   \hbox{\bf \messagetext}}
\def\filxymatrix#1{\myxymessage}
\def\filxyarray#1{\myxymessage}

\newdir^{ (}{{}*!/-3pt/\dir^{(}}
\newdir_{ (}{{}*!/-3pt/\dir_{(}}
\newdir^{ )}{{}*!/+3pt/\dir^{)}}
\newdir_{ )}{{}*!/+3pt/\dir_{)}}

\def\rscript#1{\hbox to 0pt{$\scriptstyle#1$\hss}}

\let\oldbullet\bullet
\def\bullet{{\mathchoice{\oldbullet}%
                        {\oldbullet}%
                        {\scriptscriptstyle\oldbullet}%
                        {\oldbullet}}}

\newcommand{\argdot}{{\;\bullet\;}}

\newcommand{\iat}{\mathop{i_\tau}\nolimits}

\newcommand{\CHFYAbM}{\mathop{CHM^{Ab}(Y)_F}\nolimits}
\newcommand{\CHFYPAbM}{\mathop{CHM^{Ab}_\Phi(Y)_F}\nolimits}
\newcommand{\CHFYpAbM}{\mathop{CHM^{Ab}(Y_\varphi)_F}\nolimits}

\newcommand{\CHTFSSM}{\mathop{CHMT_{\FS}(S(\FS))_F}\nolimits}
\newcommand{\CHFYpMs}{\mathop{CHM^s(\Yp)_F}\nolimits}
\newcommand{\DBcM}{\mathop{DM_{\text{\cyrb},c}}\nolimits}
\newcommand{\DBcFbM}{\mathop{\DBcM(\bullet)_F}\nolimits}

\newcommand{\DBcFBsM}{\mathop{\DBcM(\Bs)_F}\nolimits}

\newcommand{\DBcFSSM}{\mathop{\DBcM(S(\FS))_F}\nolimits}
\newcommand{\DBcXM}{\mathop{\DBcM(X)}\nolimits}
\newcommand{\DBcFXM}{\mathop{\DBcM(X)_F}\nolimits}
\newcommand{\DBcFYM}{\mathop{\DBcM(Y)_F}\nolimits}
\newcommand{\DBcFYAbM}{\mathop{DM_{\text{\cyrb},c}^{Ab}(Y)_F}\nolimits}
\newcommand{\DBcFYPAbM}{\mathop{DM_{\text{\cyrb},c,\Phi}^{Ab}(Y)_F}\nolimits}
\newcommand{\DBcFYPkAbM}{\mathop{DM_{\text{\cyrb},c,\Phi}^{Ab}(Y_k)_F}\nolimits}
\newcommand{\DBcFYpAbM}{\mathop{DM_{\text{\cyrb},c}^{Ab}(Y_\varphi)_F}\nolimits}
\newcommand{\DBcFYppAbM}{\mathop{DM_{\text{\cyrb},c}^{Ab}(Y_\varphi)_F'}\nolimits}

\newcommand{\DBcZM}{\mathop{\DBcM(Z)}\nolimits}

\newcommand{\DFTBsM}{\mathop{DMT(\Bs)_F}\nolimits}
\newcommand{\DFTSsM}{\mathop{DMT(\Ss)_F}\nolimits}
\newcommand{\DFTSSM}{\mathop{DMT_{\FS}(S(\FS))_F}\nolimits}
\newcommand{\DFTSSfM}{\mathop{DMT_{\FS}(S(\FS))_F^\natural}\nolimits}
\newcommand{\DBcFYPM}{\mathop{\DBcM(Y(\Phi))_F}\nolimits}

\newcommand{\MHM}{\mathop{\bf MHM}\nolimits}

\newcommand{\Bs}{\mathop{B_\sigma}\nolimits}
\newcommand{\SH}{\mathop{{\bf SH}}\nolimits}
\newcommand{\Ss}{\mathop{S_\sigma}\nolimits}
\newcommand{\bSs}{\mathop{\overline{\Ss}}\nolimits}
\newcommand{\St}{\mathop{S_\tau}\nolimits}
\newcommand{\Yp}{\mathop{Y_\varphi}\nolimits}
\newcommand{\pis}{\mathop{\pi_\sigma}\nolimits}
\newcommand{\pios}{\mathop{\pi'_\sigma}\nolimits}
\newcommand{\pits}{\mathop{\pi''_\sigma}\nolimits}

\newcommand{\one}{\mathds{1}}

\begin{document}

%

\hfuzz=3pt
\overfullrule=10pt                   


\setlength{\abovedisplayskip}{6.0pt plus 3.0pt}
\setlength{\belowdisplayskip}{6.0pt plus 3.0pt}
\setlength{\abovedisplayshortskip}{6.0pt plus 3.0pt}
\setlength{\belowdisplayshortskip}{6.0pt plus 3.0pt}

\setlength{\baselineskip}{13.0pt}
\setlength{\lineskip}{0.0pt}
\setlength{\lineskiplimit}{0.0pt}

%
%

\title{Weights and conservativity
\forget{
\footnotemark
\footnotetext{To appear in ....}
}
}
\author{\footnotesize by\\ \\
\mbox{\hskip-2cm
\begin{minipage}{6cm} \begin{center} \begin{tabular}{c}
J\"org Wildeshaus \footnote{
Partially supported by the \emph{Agence Nationale de la
Recherche}, project ``R\'egulateurs
et formules explicites''. }\\[0.2cm]
\footnotesize Universit\'e Paris 13\\[-3pt]
\footnotesize Sorbonne Paris Cit\'e \\[-3pt]
\footnotesize LAGA, CNRS (UMR~7539)\\[-3pt]
\footnotesize F-93430 Villetaneuse\\[-3pt]
\footnotesize France\\
{\footnotesize \tt wildesh@math.univ-paris13.fr}
\end{tabular} \end{center} \end{minipage}
\hskip-2cm}
\\[2.5cm]
}
\date{June 26, 2017}
\maketitle
\begin{abstract}
\noindent
The purpose of this article is to study conservativity in the context
of triangulated categories equipped with a weight structure.
As application, we establish (weight) conservativity for
the restriction of the (generic) $\ell$-adic
realization to the category of motives of Abelian type
of characteristic zero. \\


\end{abstract}


\bigskip
\bigskip
\bigskip

\noindent {\footnotesize Math.\ Subj.\ Class.\ (2010) numbers: 
14F20 (14C15, 14C25, 14F25, 14F42, 18A22, 19E15).
}

\eject

\tableofcontents

\bigskip


%
%

\setcounter{section}{0}
\section{Introduction}
\label{Intro}



The aim of this article is to provide a proof of \emph{weight
conservativity} of the restriction of the (generic) \emph{$\ell$-adic
realization} $R_\ell$ to the category of \emph{motives of Abelian type}
of characteristic zero (Theorem~\ref{1M}).\\

Here, weight conservativity refers to the following refinement of the usual categorical notion of conservativity: the functor $R_\ell$ detects not
only isomorphisms, but also weights (and hence, their absence).
In other words, a motive $M$ of Abelian type \emph{is without weights
$\alpha,\alpha+1,\ldots,\beta$}, for integers $\alpha \le \beta$,
in the sense of \cite[Def.~1.10]{W4} if and only if the same is true
for $R_\ell(M)$. This provides the main motivation for our study:
as shown in \loccit , absence of certain weights in the \emph{boundary
motive} of a scheme $X$ allows for the construction of its \emph{interior
motive}. \\

Weight conservativity for motives of Abelian type over a point was previously
established in \cite[Sect.~1]{W8}. This special case suffices \emph{e.g.}
for the analysis of weights in the boundary motive of \emph{Picard surfaces}
\cite[Sect.~3]{W8}, essentially because the complement of the latter
in their \emph{Baily--Borel compactifications} is of dimension zero.
The study of the boundary motive of \emph{Shimura varieties}, whose
boundary is higher dimensional, requires a version of weight conservativity
for motives over higher dimensional bases, hence the need for 
Theorem~\ref{1M}. For an application to the case of \emph{Siegel threefolds},
we refer to \cite{W11}. \\

In order to prove Theorem~\ref{1M} and the intermediate conservativity results
leading up to it (Theorems~\ref{1K} and \ref{1L}), it turns out that the 
notion of \emph{weight structure} \cite{Bo} is central. In fact, the formal
structure of the proofs is best understood in that abstract setting.
This explains the title of the present work, and also, its organization.  
Section~\ref{0} is entirely situated in the context of functors $r$,
whose source is a triangulated category $\CC$ equipped with a weight
structure $w$. The following question appears natural: assuming that the weight 
structure is \emph{bounded}, and $r$ is \emph{weight exact}, does
conservativity of the restriction of $r$ to the \emph{heart} $\CC_{w = 0}$
imply (weight) conservativity of $r$~? 
The analogous question for $t$-structures can be asked; there, 
the answer is obviously positive. Weight structures being much less rigid than 
$t$-structures, it is not surprising to see that additional hypotheses
are needed. The main such condition concerns the heart: all our abstract
results (Theorems~\ref{0C}, \ref{0E}, \ref{0G} and \ref{0H})
require $\CC_{w = 0}$ to be \emph{semi-primary} \cite[D\'ef.~2.3.1]{AK}.
The reason is that all proofs make a systematic use of the 
\emph{minimal weight filtration} \cite[Sect.~2]{W9}, whose
existence is guaranteed only if $\CC_{w = 0}$ is semi-primary. \\

This explains why we restrict our attention to motives of Abelian type:
while relative Chow motives may always be expected to form
a semi-primary category \cite[Conj.~3.4]{W9}, we are at present
far from having a general proof at our disposal. Section~\ref{1}
therefore sets up a precise motivic setting in which semi-primality
is guaranteed (Theorem~\ref{1I}). \\

Section~\ref{2} then leads up to Theorem~\ref{1M}. \\ 

One word about the Hodge theoretical picture: an analogue of
our main result Theorem~\ref{1M} should
certainly hold for the (generic) \emph{Hodge theoretical
realization} $R_{\bf{H}}$. Unfortunately,
this realization is at present not fully available. Note however that 
\cite{I} provides a Hodge theoretical realization for schemes which
are smooth over $\BC$. For the applications we have in mind \cite{W11},
it will however be necessary to consider singular schemes. 
Similarly, we need the compatibility of $R_{\bf{H}}$
with the functors $f^*$, $f_*$, $f_!$, $f^!$.
In the present work, we thus replace $R_{\bf{H}}$ by the (generic)
\emph{Betti realization}, the price to be paid being that we 
\emph{a priori} lose the intrinsic notion of weights on the target
of the realization. Nonetheless, we have sufficient control on
the situ\-ation to prove the analogue of Theorem~\ref{1L}: 
as $R_\ell$, the restriction of the (generic) Betti
realization to the category of motives of Abelian type
is conservative. \\

I wish to thank F.~D\'eglise 
for useful discussions and comments, and F.~Lecomte and A.~Papadopoulos
for the invitation to the \emph{95e Rencontre entre math\'ematiciens et 
physiciens th\'eoriciens~: ``G\'eom\'etrie, arithm\'etique et physique~:
autour des motifs''} (IRMA Strasbourg, May 2015), 
where the main results of the present article were presented. 
The observations of the re\-fe\-ree helped to improve its legibility. \\

{\bf Conventions}: Throughout the article, 
$F$ denotes a finite direct product of fields
of characteristic zero.
We fix a base scheme
$\BB$, which is of finite type over some excellent scheme 
of dimension at most two. By definition, \emph{schemes} are  
$\BB$-schemes which are separated and 
of finite type (in particular, they are excellent, and
Noetherian of finite dimension), \emph{morphisms} between schemes
are separated morphisms of $\BB$-schemes, and 
a scheme is \emph{nilregular}
if the underlying reduced scheme is regular in the usual sense. \\

We use the triangulated, $\BQ$-linear categories
$\DBcXM$ of \emph{constructible Beilinson motives} over $X$ 
\cite[Def.~15.1.1]{CD},  
indexed by schemes $X$ (always in the sense of the above conventions). 
In order to have an $F$-linear theory at one's disposal, 
one re-does the construction, but using $F$
instead of $\BQ$ as coefficients \cite[Sect.~15.2.5]{CD}. This yields 
triangulated, $F$-linear categories
$\DBcFXM$ satisfying the $F$-linear analogues of the properties
of $\DBcXM$. In particular, these categories 
are pseudo-Abelian (see \cite[Sect.~2.10]{H}).
Furthermore, the canonical functor $\DBcXM \otimes_\BQ F \to \DBcFXM$
is fully faithful \cite[Sect.~14.2.20]{CD}.
As in \cite{CD}, the symbol $\one_X$
is used to denote the unit for the tensor product in $\DBcFXM$.
We shall employ the full formalism of six operations developed in
\loccit . The reader may choose to consult \cite[Sect.~2]{H} or
\cite[Sect.~1]{W10} for concise presentations of this formalism.


\bigskip

%
%

\section{Conservativity and weight conservativity}
\label{0}



We make free use of the terminology of and basic results on weight structures
\cite[Sect.~1.3]{Bo3}.  

\begin{Def} \label{0A}
Let $r: \CC_1 \to \CC_2$
be an $F$-linear exact functor between $F$-linear triangulated cate\-go\-ries
equipped with weight structures $(\CC_{1,w \le 0} , \CC_{1,w \ge 0})$
and $(\CC_{2,w \le 0} , \CC_{2,w \ge 0})$. \\[0.1cm]
(a)~The functor $r$ is said to be \emph{weight exact} if
\[
r \bigl( \CC_{1,w \le 0} \bigr) \subset \CC_{2,w \le 0} \quad \text{and} \quad
r \bigl( \CC_{1,w \ge 0} \bigr) \subset \CC_{2,w \ge 0} \; .
\]
(b)~If $r$ is weight exact, then we denote by $r_{w = 0}$ the functor
\[
\CC_{1,w = 0} \longto \CC_{2,w = 0}
\]
induced by the restriction of $r$ to the heart $\CC_{1,w = 0}$.
\end{Def}

\begin{Lem} \label{0B}
Let $r: \CC_1 \to \CC_2$
be an $F$-linear exact functor between $F$-linear triangulated cate\-go\-ries
equipped with weight structures. We assume the following.
\begin{enumerate}
\item[(1)] The functor $r$ is weight exact.
\item[(2)] The functor $r_{w = 0}: \CC_{1,w = 0} \to \CC_{2,w = 0}$ is full.
\end{enumerate}
Then for any integer $n$, and any two objects $X \in \CC_{1,w \le n}$ and
$Z \in \CC_{1,w \ge n}$, the map
\[
r: \Hom_{\CC_1} (X,Z) \longto \Hom_{\CC_2} (r(X),r(Z))
\]
is surjective.
\end{Lem}

\begin{Proof}
Fix weight filtrations
\[
X_{\le n-1} \longto X \longto X_n \longto X_{\le n-1} [1]
\]
and
\[
Z_n \longto Z \longto Z_{\ge n+1} \longto Z_n [1]
\]
of $X$ and $Z$, respectively, with $X_n, Z_n \in \CC_{1,w = n}$,
$X_{\le n-1} \in \CC_{1,w \le n-1}$, and 
$Z_{\ge n+1} \in \CC_{1,w \ge n+1}$. Given assumption~(1),
their images under $r$ are weight filtrations of $r(X)$ and $r(Z)$
of the same type. 

Let $\beta: r(X) \to r(Z)$ be a morphism. By orthogonality for the weight structure on $\CC_2$,
the morphism $\beta$ factors through a morphism $\beta' : r(X_n) \to r(Z_n)$.
The shift by $[-n]$ of the latter gives $\beta'[-n] : r(X_n)[-n] \to r(Z_n)[-n]$,
a morphism in $\CC_{2,w = 0}$. Given assumption~(2), the morphism $\beta'[-n]$
lies in the image of $r$. But then so do $\beta'$ and $\beta$.
\end{Proof}

Before stating our first abstract result on conservativity, let us recall
an important notion.

\begin{Def}
Let $\CC$ be an $F$-linear triangulated cate\-go\-ry
equipped with a weight structure $w$. Let $M \in \CC$, and $n \in \BZ$. 
A \emph{minimal weight filtration concentrated at $n$} of $M$ 
is a weight filtration
\[
M_{\le n-1} \longto M \longto M_{\ge n} \stackrel{\delta}{\longto} M_{\le n-1}[1]
\]
($M_{\le n-1} \in \CC_{w \le n-1}$, $M_{\ge n} \in \CC_{w \ge n}$) 
such that the morphism $\delta$ belongs to the \emph{radical} 
\cite[D\'ef.~1.4.1]{AK} of $\CC$: 
\[
\delta \in \rad_\CC (M_{\ge n},M_{\le n-1}[1]) \; .
\]
\end{Def}

Any two minimal weight filtrations of the same object $M$ are related
by an isomorphism (which in general is \emph{not} unique) 
\cite[proof of Thm.~2.2~(b)]{W9}.

\begin{Prop} \label{0D}
Let $r: \CC_1 \to \CC_2$
be an $F$-linear exact functor between $F$-linear triangulated cate\-go\-ries
equipped with weight structures. We assume the following.
\begin{enumerate}
\item[(1)] The functor $r$ is weight exact.
\item[(2)] The functor $r_{w = 0}: \CC_{1,w = 0} \to \CC_{2,w = 0}$ is full.
\end{enumerate}
Then $r$ maps minimal weight filtrations to minimal weight filtrations.
\end{Prop}

\begin{Proof}
Let 
\[
(\ast) \quad\quad\quad\quad
M_{\le n-1} \longto M \longto M_{\ge n} \stackrel{\delta}{\longto} M_{\le n-1}[1]
\]
be a minimal weight filtration concentrated at $n$ of $M \in \CC_1$.
Thanks to assumption~(1), the exact triangle $r(\ast)$ is a weight filtration
of $r(M)$. 

According to Lemma~\ref{0B}, any morphism 
$r(M_{\le n-1}[1]) \to r(M_{\ge n})$ lies in the image of $r$. 
Using this information, and
applying the definition, one sees
that $r(\delta)$ belongs to the radical of $\CC_2$ since $\delta$
belongs to the radical of $\CC_1$. Therefore, $r(\ast)$ is a 
minimal weight filtration concentrated at $n$ of $r(M)$. 
\end{Proof}

\begin{Thm} \label{0C}
Let $r: \CC_1 \to \CC_2$
be an $F$-linear exact functor between $F$-linear triangulated cate\-go\-ries
equipped with weight structures. We assume the following.
\begin{enumerate}
\item[(1)] The weight structure on $\CC_1$ is bounded, \emph{i.e.},
its heart $\CC_{1,w = 0}$ generates $\CC_1$ as a triangulated category.
\item[(2)] The heart $\CC_{1,w = 0}$ is \emph{semi-primary}
\cite[D\'ef.~2.3.1]{AK} and pseudo-Abelian. 
\item[(3)] The functor $r$ is weight exact.
\item[(4)] The functor $r_{w = 0}$ is full.
\item[(5)] The functor $r_{w = 0}$ is conservative.
\end{enumerate}
Then $r$ is conservative.
\end{Thm}

\begin{Proof}
The functor $r$ being an exact functor between triangulated ca\-te\-gories, it suffices to show
that only the zero object is mapped to zero under $r$. Thus, let $M \in \CC_1$,
and assume that $r(M) = 0$. 

Given assumption~(1), the category $\CC_1$ is pseudo-Abelian if and only if
its heart $\CC_{1,w = 0}$ is \cite[Lemma~5.2.1]{Bo}. Given assumption~(2),
we are thus in the abstract
situation studied in \cite[Sect.~2]{W9}, meaning in particular that 
minimal weight filtrations \emph{do} exist 
for any object of $\CC_1$ \cite[Thm.~2.2~(a)]{W9}. Let
\[
(\ast) \quad\quad\quad\quad
M_{\le -1} \longto M \longto M_{\ge 0} \stackrel{\delta}{\longto} M_{\le -1}[1]
\]
be a minimal weight filtration of $M$, \emph{i.e.},
$M_{\le -1} \in \CC_{1,w \le -1}$,
$M_{\ge 0} \in \CC_{1,w \ge 0}$,
and the morphism $\delta$ belongs to the radical of $\CC_1$.

On the one hand, 
the triangle $r(\ast)$ is exact, and
$r(M) = 0$. Therefore, the morphism $r(\delta)$ is an isomorphism.

On the other hand, according to Proposition~\ref{0D} (applicable thanks to
hypotheses~(3) and (4)), $r(\delta)$ belongs to the radical of $\CC_2$. 

But then $r(\delta)$ must be the zero morphism $0 \to 0$, meaning that
both $M_{\le -1}$ and $M_{\ge 0}$ are mapped to zero under $r$. 

According to hypothesis~(1), the weight
structure on $\CC_1$ is bounded.
Therefore, the above procedure, successively applied to minimal weight
filtrations concentrated at integers different from zero, allows to
reduce to the case where $M$ is pure of some weight, say $n$. But then
$M[-n]$ belongs to $\CC_{1,w = 0}$. Now apply assumption~(5).
\end{Proof}

\begin{Rem}
The analogue of Theorem~\ref{0C} for triangulated categories equipped
with $t$-structures holds, and requires only the analogues of assumptions~(1),
(3), and (5).
\end{Rem}

Recall the following.

\begin{Def}[{\cite[Def.~1.6, Def.~1.10]{W4}}] 
Let $\CC$ be an $F$-linear triangulated cate\-go\-ry
equipped with a weight structure $w$. Let $M \in \CC$, and $\alpha \le \beta$
two integers (which may be identical). \\[0.1cm]
(a)~A \emph{weight filtration of $M$ avoiding weights 
$\alpha,\alpha+1,\ldots,\beta$} is an exact triangle
\[ 
M_{\le {\alpha-1}} \longto M \longto M_{\ge {\beta+1}} \longto M_{\le \alpha-1}[1]
\]
in $\CC$, with $M_{\le {\alpha-1}} \in \CC_{w \le \alpha-1}$
and $M_{\ge {\beta+1}} \in \CC_{w \ge {\beta+1}} \ $. \\[0.1cm]
(b)~We say that
$M \in \CC$ \emph{does not have weights $\alpha,\alpha+1,\ldots,\beta$},
or that $M$ \emph{is without weights $\alpha,\alpha+1,\ldots,\beta$},
if it admits a weight filtration 
avoiding weights $\alpha,\alpha+1,\ldots,\beta$.
\end{Def}

We leave it to the reader to verify that any weight filtration avoiding weights 
$\alpha,\alpha+1,\ldots,\beta$ is a minimal weight filtration concentrated at $n$,
for any integer $n$ between $\alpha$ and $\beta+1$. \\

The following result sharpens Theorem~\ref{0C}; it will
however not be needed in the sequel.

\begin{Thm} \label{0E}
Let $r: \CC_1 \to \CC_2$
be an $F$-linear exact functor between $F$-linear triangulated cate\-go\-ries
equipped with weight structures. The assumptions (1)--(5) are the same as
in Theorem~\ref{0C}.
Then $r$ is \emph{weight conservative}, \emph{i.e.}, it detects weights.
More precisely, let $M \in \CC_1$, and $\alpha \le \beta$ two integers. \\[0.1cm]
(a)~$M$ lies in the heart $\CC_{1,w = 0}$ if and only if 
$r(M)$ lies in the heart $\CC_{2,w = 0}$. \\[0.1cm]
(b)~$M$ lies in $\CC_{1,w \le \alpha}$ if and only if 
$r(M)$ lies in $\CC_{2,w \le \alpha}$. \\[0.1cm] 
(c)~$M$ lies in $\CC_{1,w \ge \beta}$ if and only if 
$r(M)$ lies in $\CC_{2,w \ge \beta}$.. \\[0.1cm]
(d)~$M$ is without weights $\alpha,\alpha+1,\ldots,\beta$ if and only if 
$r(M)$ is without weights $\alpha,\alpha+1,\ldots,\beta$. 
\end{Thm}

\begin{Proof}
According to assumption~(1), the weight structure 
on $\CC_1$ is boun\-ded, and the ``only if'' parts of
statements~(a)--(d) are true. It therefore suffices to prove the
``if'' part of statement~(d).

Consider minimal weight filtrations concentrated at weight
$\alpha$ and at weight $\beta+1$, respectively: 
\[
M_{\le \alpha-1} \longto M \longto M_{\ge \alpha} 
\stackrel{\delta_\alpha}{\longto} M_{\le \alpha-1}[1] \; ,
\]
\[
M_{\le \beta} \longto M \longto M_{\ge \beta+1} 
\stackrel{\delta_{\beta+1}}{\longto} M_{\le \beta}[1] \; .
\]
By orthogonality, the identity on $M$ extends to a morphism of exact
triangles
\[
\vcenter{\xymatrix@R-10pt{
        M_{\le \alpha-1} \ar[d]_{m} \ar[r] & M \ar@{=}[d] \ar[r] &
        M_{\ge \alpha} \ar[d] \ar[r]^-{\delta_\alpha} & M_{\le \alpha-1}[1] \ar[d]_{m[1]}\\
        M_{\le \beta} \ar[r] & M \ar[r] &
        M_{\ge \beta+1} \ar[r]^-{\delta_{\beta+1}} & M_{\le \beta}[1]
\\}}
\] 
By Proposition~\ref{0D}, the image under $r$ of the above is a morphism relating 
minimal weight filtrations concentrated at weight
$\alpha$ and at weight $\beta+1$, respectively:
\[
\vcenter{\xymatrix@R-10pt{
        r \bigl( M_{\le \alpha-1} \bigr) \ar[d]_{r(m)} \ar[r] & 
        r(M) \ar@{=}[d] \ar[r] &
        r \bigl( M_{\ge \alpha} \bigr)\ar[d] \ar[r]^-{r(\delta_\alpha)} & 
        r \bigl( M_{\le \alpha-1}[1]  \bigr) \ar[d]_{r(m)[1]}\\
        r \bigl( M_{\le \beta} \bigr) \ar[r] & 
        r(M) \ar[r] &
        r \bigl( M_{\ge \beta+1} \bigr) \ar[r]^-{r(\delta_{\beta+1})} & 
        r \bigl( M_{\le \beta}[1] \bigr)
\\}}
\] 
Thus, both $r(\delta_\alpha)$ and $r(\delta_{\beta+1})$ 
lie in the radical of $\CC_2$.

But $r(M)$ is suppposed to be without weights $\alpha,\alpha+1,\ldots,\beta$,
meaning that there is a minimal weight filtration of $r(M)$, which is 
concentrated at any integer between $\alpha$ and $\beta+1$.
Unicity of minimal weight filtrations in $\CC_2$
\cite[proof of Thm.~2.2~(b)]{W9} shows that this latter weight filtration
is isomorphic to both the top and the bottom of the above diagram, meaning
in particular that they are (abstractly) isomorphic to each other.
In particular, the object $r(M_{\le \beta})$ belongs to $\CC_{2,w \le \alpha-1}$.
Orthogonality then allows to extend the identity on $r(M)$ to a morphism of exact
triangles
\[
\vcenter{\xymatrix@R-10pt{
        r \bigl( M_{\le \beta} \bigr) \ar[d]_{n} \ar[r] & 
        r(M) \ar@{=}[d] \ar[r] &
        r \bigl( M_{\ge \beta+1} \bigr) \ar[d] \ar[r]^-{r(\delta_{\beta+1})} & 
        r \bigl( M_{\le \beta}[1] \bigr) \ar[d]_{n[1]} \\
        r \bigl( M_{\le \alpha-1} \bigr) \ar[r] & 
        r(M) \ar[r] &
        r \bigl( M_{\ge \alpha} \bigr) \ar[r]^-{r(\delta_\alpha)} & 
        r \bigl( M_{\le \alpha-1}[1]  \bigr)\\
\\}}
\] 
Using the fact that both $r(\delta_\alpha)$ and $r(\delta_{\beta+1})$ 
lie in the radical, one shows that both compositions $n \circ r(m)$ and
$r(m) \circ n$ are automorphisms; in particular, $r(m)$ is an isomorphism. 

But then (Theorem~\ref{0C}), so is $m$ itself. This yields a weight filtration
\[
M_{\le \alpha-1} \longto M \longto M_{\ge \beta+1} 
\longto M_{\le \alpha-1}[1] 
\]
of $M$ avoiding weights $\alpha,\alpha+1,\ldots,\beta$.
\end{Proof}

The second half of the present section treats conservativity
in a slightly different context. The source of the functors 
in question remains equipped with a weight structure. But their target
is only supposed to be Abelian.
Recall \cite[Prop.~2.1.2~1]{Bo} that any (covariant) additive functor
$\CH$ from a triangulated category $\CC$ carrying a weight structure $w$,
to an Abelian category $\FA$ admits a canonical \emph{weight filtration}
by sub-functors
\[
\ldots \subset W_n \CH \subset W_{n+1} \CH \subset \ldots \subset \CH \; .
\]
According to \cite[Def.~2.1.1]{Bo} (use the normalization of \cite[Def.~1.3.1]{Bo3} for the
signs of the weights), for an object $M$ of $\CC$, and $n \in \BZ$,
the sub-object $W_n \CH(M) \subset \CH(M)$ is defined as the image of the morphism
$\CH(\iota_{w \le n})$, for \emph{any} weight filtration
\[
M_{w \le n} \stackrel{\iota_{w \le n}}{\longto} M \longto M_{w \ge n+1} \longto M_{w \le n}[1]
\]
(with $M_{w \le n} \in \CC_{w \le n}$ and $M_{w \ge n+1} \in \CC_{w \ge n+1}$).  
For any $m \in \BZ$, one defines 
\[
\CH^m : \CC \longto \FA \; , \; M \longto \CH \bigl( X[m] \bigr) \; ;
\]
according to the usual convention, the weight filtration of $\CH^m(M)$ \emph{equals} the
weight filtration of $\CH(X[m])$, \emph{i.e.}, it differs by \emph{d\'ecalage}
from the intrinsic weight filtration of the covariant additive functor $\CH^m$. \\

The proof of the following result is formal, and therefore left to the reader
(cmp.\ \cite[Lemma~1.11]{W8}).

\begin{Lem} \label{0F}
Let $\CH: \CC_1 \to \CC_3$
be an $F$-linear functor 
from an $F$-linear triangulated category $\CC_1$ equipped with a weight structure
to an $F$-linear Abelian cate\-go\-ry $\CC_3$. We assume the following.
\begin{enumerate}
\item[(1)] Any morphism in the image of $\CH$ is strict with respect to the
weight filtration of $\CH$.
\item[(2)] The restriction of $\CH$ to the heart
$\CC_{1,w = 0}$ maps the radical to zero.
\end{enumerate}
Then for any integer $n$, and any two objects $X \in \CC_{1,w \le n}$ and
$Z \in \CC_{1,w \ge n}$, the map
\[
\CH: \Hom_{\CC_1} (Z,X) \longto \Hom_{\CC_3} \bigl( \CH(Z),\CH(X) \bigr)
\]
maps the radical to zero.
\end{Lem}

\begin{Thm} \label{0G}
Let $\CH: \CC_1 \to \CC_3$
be an $F$-linear homological functor 
from an $F$-linear triangulated category $\CC_1$ equipped with a weight structure
to an $F$-linear Abelian cate\-go\-ry $\CC_3$. We assume the following.
\begin{enumerate}
\item[(1)] The weight structure on $\CC_1$ is bounded.
\item[(2)] The heart $\CC_{1,w = 0}$ is semi-primary and pseudo-Abelian. 
\item[(3)] Any morphism in the image of $\CH$ is strict with respect to the
weight filtration of $\CH$.
\item[(4)] The restriction of $\CH$ to the heart
$\CC_{1,w = 0}$ maps the radical to zero.
\item[(5)] Zero is the only object of the heart $\CC_{1,w = 0}$
mapped to zero under (the restriction to $\CC_{1,w = 0}$ of) 
the functor $(\CH^m)_{m \in \BZ}$.
\end{enumerate}
Then $(\CH^m)_{m \in \BZ}$ is conservative.
\end{Thm}

\begin{Proof}
Let $M \in \CC_1$. Let
\[
(\ast) \quad\quad\quad\quad
M_{\le -1} \longto M \longto M_{\ge 0} \stackrel{\delta}{\longto} M_{\le -1}[1]
\]
be a minimal weight filtration of $M$ (here, we use assumptions~(1) and (2)).
According to Lemma~\ref{0F} (applicable thanks to
hypotheses~(3) and (4)), $\CH(\delta) = 0$. The functor $\CH$ is supposed
to be homological, hence if $\CH(M)$ is zero, then
both $\CH(M_{\le -1})$ and $\CH(M_{\ge 0})$ are zero. 
Similarly, if $\CH^m(M)$ is zero for some $m \in \BZ$, then
so are $\CH^m(M_{\le -1})$ and $\CH^m(M_{\ge 0})$.

Now let $M \in \CC_1$, and assume that $\CH^m(M)=0$ for all $m \in \BZ$.
We need to show that $M = 0$.
According to hypothesis~(1), the weight
structure on $\CC_1$ is bounded.
Therefore, the above procedure, successively applied to minimal weight
filtrations concentrated at integers different from zero, allows to
reduce to the case where $M$ is pure of some weight, say $n$. But then
$M[-n]$ belongs to $\CC_{1,w = 0}$. Now apply assumption~(5).
\end{Proof}

As before, there is a version ``with weights'' of the above.

\begin{Thm} \label{0H}
Let $\CH: \CC_1 \to \CC_3$
be an $F$-linear homological functor 
from an $F$-linear triangulated category $\CC_1$ equipped with a weight structure
to an $F$-linear Abelian cate\-go\-ry $\CC_3$. The assumptions (1)--(5) are the same as
in Theorem~\ref{0G}.
Then $(\CH^m)_{m \in \BZ}$ is weight conservative, \emph{i.e.}, it detects weights.
More precisely, let $M \in \CC_1$, and $\alpha \le \beta$ two integers. \\[0.1cm]
(a)~$M$ lies in the heart $\CC_{1,w = 0}$ if and only if 
$\CH^n (M)$ is pure of weight $n$, for all $n \in \BZ$. \\[0.1cm]
(b)~$M$ lies in $\CC_{1,w \le \alpha}$ if and only if 
$\CH^n (M)$ is of weights $\le n + \alpha$, for all $n \in \BZ$. \\[0.1cm] 
(c)~$M$ lies in $\CC_{1,w \ge \beta}$ if and only if 
$\CH^n (M)$ is of weights $\ge n + \beta$, for all $n \in \BZ$. \\[0.1cm]
(d)~$M$ is without weights $\alpha,\alpha+1,\ldots,\beta$ if and only if 
$\CH^n (M)$ is without weights $n + \alpha, n + \alpha+1,\ldots, n + \beta$, 
for all $n \in \BZ$.  
\end{Thm}

\begin{Proof}
According to assumption~(1), the weight structure 
on $\CC_1$ is boun\-ded, and the ``only if'' parts of
statements~(a)--(d) are true. It therefore suffices to prove the
``if'' part of statement~(d).

Consider minimal weight filtrations concentrated at weight
$\alpha$ and at weight $\beta+1$, respectively: 
\[
M_{\le \alpha-1} \longto M \longto M_{\ge \alpha} 
\stackrel{\delta_\alpha}{\longto} M_{\le \alpha-1}[1] \; ,
\]
\[
M_{\le \beta} \longto M \longto M_{\ge \beta+1} 
\stackrel{\delta_{\beta+1}}{\longto} M_{\le \beta}[1] \; .
\]
By orthogonality, the identity on $M$ extends to a morphism of exact
triangles
\[
\vcenter{\xymatrix@R-10pt{
        M_{\le \alpha-1} \ar[d]_{m} \ar[r] & M \ar@{=}[d] \ar[r] &
        M_{\ge \alpha} \ar[d] \ar[r]^-{\delta_\alpha} & M_{\le \alpha-1}[1] \ar[d]_{m[1]}\\
        M_{\le \beta} \ar[r] & M \ar[r] &
        M_{\ge \beta+1} \ar[r]^-{\delta_{\beta+1}} & M_{\le \beta}[1]
\\}}
\] 
By Lemma~\ref{0F}, both $\CH^* (\delta_\alpha)$ and 
$\CH^* (\delta_{\beta+1})$ are zero. Thus, the above morphism of exact
triangles induces a morphism of exact sequences
\[
\vcenter{\xymatrix@R-10pt{
        0 \ar[r] &
        \CH^* (M_{\le \alpha-1}) \ar@{^{ (}->}[d]_{\CH^*(m)} \ar[r] & 
        \CH^* (M) \ar@{=}[d] \ar[r] &
        \CH^* (M_{\ge \alpha}) \ar@{->>}[d] \ar[r] & 0 \\
        0 \ar[r] &
        \CH^* (M_{\le \beta}) \ar[r] & 
        \CH^* (M) \ar[r] &
        \CH^* (M_{\ge \beta+1}) \ar[r] & 0
\\}}
\] 
Our hypothesis on weights avoided in $\CH^* (M)$ implies that
the monomorphism $\CH^* (m)$ is in fact an isomorphism.

But then (Theorem~\ref{0G}), so is $m$ itself. This yields a weight filtration
\[
M_{\le \alpha-1} \longto M \longto M_{\ge \beta+1} 
\longto M_{\le \alpha-1}[1] 
\]
of $M$ avoiding weights $\alpha,\alpha+1,\ldots,\beta$.
\end{Proof}


\bigskip

%
%

\section{Relative motives of Abelian type}
\label{1}



Let $S$ be a scheme (in the sense of the conventions fixed in our Introduction).

\begin{Def} \label{1A}
A \emph{good stratification} of $S$ 
indexed by a finite set $\FS$
is a collection of locally closed sub-schemes $\Ss$ indexed by
$\sigma \in \FS$, such that 
\[
S = \coprod_{\sigma \in \FS} \Ss
\]
on the set-theoretic level, and such that the closure $\bSs$
of any stratum $\Ss$ is a union of strata $S_\tau \, $. 
\end{Def}

In the setting of Definition~\ref{1A},
we shall often write $S(\FS)$ instead of $S$. Recall the following result.

\begin{Thm}[{\cite[Thm.~4.5~(a)]{W9}}] \label{1B}
Let $S(\FS) = \coprod_{\sigma \in \FS} \Ss$ be a good stratification
of a scheme $S(\FS)$. Assume the following for all $\sigma \in \FS$:
($\alpha$)~the stratum $\Ss$ is nilregular, ($\beta$)~for any stratum 
$\iat: \St \into \bSs$ contained in the closure $\bSs$ of $\Ss$, 
the functor
\[
i_\tau^! : \DBcM(\bSs)_F \longto \DBcM(S_\tau)_F
\]
maps $\one_{\bSs}$ to a Tate motive over $S_\tau$. Then
the categories $\DFTSsM$ of \emph{Tate motives over $\Ss$}
\cite[Sect.~3.3]{L3}, $\sigma \in \FS$, 
can be glued to give a full,
triangulated sub-category $\DFTSSM$ of $\DBcFSSM \, $.
\end{Thm}

Recall \cite[Rem.~4.7]{W9} that thanks to \emph{absolute purity}
\cite[Thm.~14.4.1]{CD}, hypotheses ($\alpha$) and ($\beta$) from Theorem~\ref{1B}
are satisfied as soon as the closures $\bSs$ of all strata
$\Ss$, $\sigma \in \FS$ are nilregular.

\begin{Def} \label{1C}
Let $S(\FS) = \coprod_{\sigma \in \FS} \Ss$ be a good stratification
of a scheme $S(\FS)$. Assume the following for all $\sigma \in \FS$:
($\alpha$)~$\Ss$ is nilregular, ($\beta$)~for any 
$\iat: \St \into \bSs$, the functor $i_\tau^!$
maps $\one_{\bSs}$ to a Tate motive over $S_\tau$. 
The category $\DFTSSM$ of Theorem~\ref{1B}
is called the category of \emph{$\FS$-constructible
Tate motives over $S(\FS)$}. 
\end{Def}

According to \cite[Thm.~4.5~(d)]{W9}, the category $\DFTSSM$
is pseudo-Abelian.
Now let $S(\FS) = \coprod_{\sigma \in \FS} \Ss$
and $Y(\Phi) = \coprod_{\varphi \in \Phi} \Yp$
be good stratifications of schemes $S(\FS)$ and $Y(\Phi)$,
respectively.

\begin{Def} \label{1D}
A morphism $\pi: S(\FS) \to Y(\Phi)$ is said to be a
\emph{morphism of good stratifications} if
the pre-image $\pi^{-1}(\Yp)$ of 
any stratum $\Yp$ of $Y(\Phi)$, $\varphi \in \Phi \, $ 
is a union of strata $\Ss \, $.
\end{Def}

\forget{
\begin{Def} \label{1E}
Let $\pi: S(\FS) \to Y(\Phi)$ be a morphism of good stratifications.
Assume the following for all $\sigma \in \FS$:
($\alpha$)~$\Ss$ is nilregular, ($\beta$)~for any 
$\iat: \St \into \bSs$, the functor $i_\tau^!$
maps $\one_{\bSs}$ to a Tate motive over $S_\tau$.
Define the category $\pi_* \! \DFTSSfM$
as the strict, full, dense, $F$-linear
triangulated sub-category of $\DBcFYPM$ generated by
the images under $\pi_*$ of the objects of $\DFTSSM \, $. 
\end{Def}
}
\begin{Def} \label{1F}
(a)~A morphism $\pi: S(\FS) \to Y(\Phi)$ of good stratifications
is said to be \emph{of Abelian type} if it is proper, and if
the following conditions are satisfied. 
\begin{enumerate}
\item[(1)] All strata $\Yp$ and $\Ss$, $\varphi \in \Phi$, $\sigma \in \FS$,
are nilregular, and for any 
$\iat: \St \into \bSs$, the functor $i_\tau^!$
maps $\one_{\bSs}$ to a Tate motive over $S_\tau$. 
\item[(2)] For all $\sigma \in \FS$
such that $\Ss$ is a stratum
of $\pi^{-1}(\Yp)$, the morphism
$\pis : \Ss \to \Yp$ can be factorized,
\[
\pis = \pios \circ \pits : \Ss \stackrel{\pits}{\longto} \Bs 
\stackrel{\pios}{\longto} \Yp \; ,
\]
such that the motive 
\[
\pi''_{\sigma,*} \one_{\Ss} \in \DBcFBsM
\]
belongs to the category $\DFTBsM$ of Tate motives over $\Bs \, $,
the morphism $\pios$ is proper and smooth, 
and its pull-back to any geometric
point of $\Yp$ lying over a generic point
is isomorphic to a finite disjoint union of Abelian varieties.
\end{enumerate}  
(b)~Let $Y$ be a scheme, equipped with a good stratification $\Phi$
with nilregular strata. 
Define the category $\DBcFYPAbM$
as the strict, full, dense, $F$-linear
triangulated sub-category of $\DBcFYM$ generated by
the images under $\pi_*$ of the objects of $\DFTSSM \, $,
where $\pi: S(\FS) \to Y(\Phi)$ runs through the morphisms of Abelian type with target
equal to $Y = Y(\Phi)$. Objects of $\DBcFYPAbM$ are called
\emph{$\Phi$-constructible motives of Abelian type over $Y$}. \\[0.1cm]
(c)~Let $Y$ be a nilregular scheme. Set
\[
\DBcFYAbM := \DBcFYPAbM \; ,
\]
for the trivial stratification $\Phi = \{ \varphi \}$.
\end{Def}

Since the category $\DBcFYM$ is pseudo-Abelian (see \cite[Sect.~2.10]{H}),
so is $\DBcFYPAbM$. \\

If a nilregular scheme $Y$ is equipped with a good stratification $\Phi$
with nilregular strata, 
then by \cite[Thm.~4.5~(b)]{W9} and \cite[Thm.~15.2.4]{CD}, the category
$\DBcFYPAbM$ is closed under duality. \\

Next, we need to discuss weight structures.

\begin{Thm} \label{1G}
Let $Y$ be a scheme, equipped with a good stratification $\Phi$
with nilregular strata. \\[0.1cm]
(a)~The \emph{motivic weight structure} $w$ on $\DBcFYPM$ 
(\cite[Thm.~3.3]{H}, \cite[Thm.~2.1.1]{Bo4}) induces a weight structure,
still denoted by the same letter $w$, on
$\DBcFYPAbM$. \\[0.1cm]
(b)~The motivic weight structure on $\DBcFYPAbM$ is bounded.
\end{Thm}

\begin{Proof}
We imitate the proof of \cite[Cor.~4.11]{W9}.
Let $\CK$ be the strict, full, $F$-linear sub-category of $\DBcFYPM$ 
of finite direct sums of motives isomorphic to images  under $\pi_*$ of objects in
\[
\CHTFSSM := \DFTSSM \cap \, CHM(S(\FS))_F \; ,
\]
for morphisms $\pi: S(\FS) \to Y(\Phi)$ of Abelian type.
Denote by $\CD$ the triangulated category generated by $\CK$.

According to our definition and \cite[Cor.~4.12~(b)]{W9}, 
$\CD$ equals the strict, full, $F$-linear
triangulated sub-category of $\DBcFYM$ generated by
the images under $\pi_*$ of the objects of $\DFTSSM \, $,
where $\pi: S(\FS) \to Y(\Phi)$ runs through the morphisms of Abelian type.

Following \cite[Rem.~4.4]{W9}, 
the motivic weight structure induces a bounded weight structure on $\CD$,
whose heart contains $\CK$. 

Repeat the same argument with the pseudo-Abelian completion
$\CK^\natural$ instead of $\CK$. We get a bounded weight structure,
induced by the motivic weight structure,
on a triangulated sub-category $\CD^\natural$ of $\DBcFYPM$.

Our claim is implied by \cite[Prop.~5.2.2]{Bo}, which states that 
$\CD^\natural$ is the pseudo-Abelian competion of $\CD$, hence equal
to $\DBcFYPAbM$.
\end{Proof}

\begin{Def} \label{1H}
(a)~Let $Y$ be a scheme, equipped with a good stratification $\Phi$
with nilregular strata.
A \emph{$\Phi$-constructible Chow motive of Abelian type over $Y$} is an object of
\[
\CHFYPAbM := DM_{\text{\cyrb},c,\Phi}^{Ab}(Y)_{F,w= 0} \; .
\]
(b)~Let $Y$ be a nilregular scheme. Set
\[
\CHFYAbM := \CHFYPAbM \; ,
\]
for the trivial stratification $\Phi = \{ \varphi \}$.
\end{Def}

Note that since $\DBcFYPAbM$ is pseudo-Abelian, so is $\CHFYPAbM$.
Using \cite[Thm.~4.3.2~II]{Bo},
let us extract the following from the proof of Theorem~\ref{1G}.

\begin{Lem} \label{1lem}
Let $Y$ be a scheme, equipped with a good stratification $\Phi$
with nilregular strata. Then the strict, full, pseudo-Abelian, $F$-linear  
sub-category $\CHFYPAbM$ of $\DBcFYM$ is generated by the
images  under $\pi_*$ of objects in $\CHTFSSM \,$,
where $\pi: S(\FS) \to Y(\Phi)$ runs through the morphisms of Abelian type
with target $Y(\Phi)$.
\end{Lem}

In order to have the main result from \cite{W9} at our disposal,
let us check its hypotheses.

\begin{Prop} \label{1Ipre}
Let $\pi: S(\FS) \to Y(\Phi)$ be a morphism of Abelian type. \\[0.1cm]
(a)~Let $\sigma \in \FS$
such that $\Ss$ is a stratum
of $\pi^{-1}(\Yp)$, and
\[
\pis = \pios \circ \pits : \Ss \stackrel{\pits}{\longto} \Bs 
\stackrel{\pios}{\longto} \Yp 
\]
a factorization of
the morphism $\pis : \Ss \to \Yp$ 
as in Definition~\ref{1F}~(a)~(2).
Then the \emph{smooth Chow motive over $\Yp$} 
\cite[Def.~5.16]{L2}
\[
\pi'_{\sigma,*} \one_{\Bs} \in \CHFYpMs 
\]
is \emph{finite dimensional} \cite[D\'ef.~9.1.1]{AK} 
(cmp.\ \cite[Def.~3.7]{Ki}). \\[0.1cm]
(b)~The morphism $\pi$ satisfies the
assumptions of \cite[Main Thm.~5.4]{W9}.
\end{Prop}

\begin{Proof}
The morphism $\pios$ is proper and smooth, 
and its pull-back to any geometric
point of $\Yp$ lying over a generic point
is isomorphic to a finite disjoint union of Abelian varieties.
According to \cite[pp.~54--55]{O'S}, finite dimensionality can be checked
after base change to the geometric generic points of $\Yp$.
Now apply \cite[Thm.~(3.3.1)]{Ku}. This establishes part~(a) of the claim.

But given Definition~\ref{1F}, part~(a) is all that is needed in order
to show that the assumptions of \cite[Main Thm.~5.4]{W9} are satisfied.
\end{Proof}

We thus get the following structural result.

\begin{Thm} \label{1I}
Let $Y$ be a scheme, equipped with a good stratification $Y = Y(\Phi)$
with nilregular strata. Then 
the $F$-category $\CHFYPAbM$ 
of $\Phi$-constructible Chow motives of Abelian type over $Y$ is semi-primary
(and pseudo-Abelian). 
\end{Thm}  

\begin{Proof}
Given \cite[Prop.~2.3.4~c)]{AK}, our claim follows from Lemma~\ref{1lem},
Proposition~\ref{1Ipre}~(b) and \cite[Main Thm.~5.4]{W9}.
\end{Proof}

Here is our first result on conservativity in the motivic context.

\begin{Thm} \label{1K}
Fix a generic point $\Spec k \into \BB$ of the base scheme $\BB$.
Let $Y = Y(\Phi) = \coprod_{\varphi \in \Phi} \Yp$ be a good stratification
with nilregular strata,
such that the generic points of all $\Yp$ lie over $\Spec k \into \BB$. 
Denote by $Y_k$ the base change of $Y$ to $\Spec k$.
Then the inverse image functor
\[
\DBcFYPAbM \longto \DBcFYPkAbM 
\]
is conservative.  
\end{Thm}

\begin{Proof}
According to \cite[Thm.~14.3.3]{CD}, the categories $\DBcM (\argdot)_F$ are 
\emph{separated} in the sense of \cite[Def.~2.1.7]{CD}.
Therefore, it suffices to check the claim after application
of the inverse image functors to all $\Yp \, $. Given proper base change 
\cite[Thm.~2.4.50~(4)]{CD}, we are thus reduced to the case where
the stratification $\Phi$ consists of a single stratum: $Y(\Phi) = \Yp \, $.
Recall that by assumption, the scheme $\Yp$ is nilregular.

According to Lemma~\ref{1lem} and \cite[Prop.~5.5]{W9},
every object of $\CHFYPAbM$
is a direct factor of a finite direct sum of objects isomorphic to
$\pi'_* \one_B(p)[2p]$, for $p \in \BZ$ and proper and smooth morphisms
$\pi' : B \to \Yp$.

In particular, $\CHFYPAbM$
is contained in the category of smooth Chow motives over $\Yp \, $.
By assumption, the morphism $Y_{\varphi,k} \to \Yp$
is dominant. Conservativity of the restriction of the inverse
image
\[
\alpha^*: \DBcFYPAbM \longto \DBcFYPkAbM
\]
to $\CHFYPAbM$ thus follows from
\cite[pp.~54--55]{O'S}.

To treat the full triangulated category $\DBcFYPAbM \,$, note that
according to Theorem~\ref{1G}~(b), its weight structure is bounded.

By Theorem~\ref{1I}, the heart of the weight structure
is semi-primary and pseudo-Abelian. 

Furthermore, the functor $\alpha^*$ is weight exact. 

According to \cite[Prop.~5.1.1]{O'S},
the restriction of the
functor $\alpha^*$ to the heart $\CHFYPAbM$
is full. 

Thus, the assumptions of Theorem~\ref{0C} are all satisfied.
\end{Proof}


\bigskip

%
%

\section{Realizations}
\label{2}



This section will be devoted to realizations. 
We fix a generic point $\Spec k$ of our base scheme $\BB$,
and assume that we are in one of the following situations.
\begin{enumerate}
\item[(i)] $k$ is embedded into $\BC$ \emph{via} a morphism 
$\eta: k \into \BC$, yielding a geometric point
of $\BB$, denoted by the same symbol
\[
\eta: \Spec \BC \longto \Spec k 
\longinto \BB \; . 
\]
The \emph{Betti realization} is defined in \cite[D\'ef.~2.1]{Ay2}. 
It is a family of covariant exact functors
\[
R_{\eta,Z}: \SH(Z) \longto D(Z) \; ,
\] 
indexed by quasi-projective $k$-schemes $Z$. 
The source of $R_{\eta,Z}$ is
the \emph{stable homotopy category of $Z$-schemes} \cite[Sect.~4.5]{Ay1}.
Its target is 
the derived category of the Abelian category 
of sheaves with values in Abelian groups
on the topological space $Z(\BC)$ of points 
of $Z$ with values in $\BC$ with respect to $\eta$. The functors $R_{\eta,Z}$ 
are symmetric monoidal \cite[Lemme~2.2]{Ay2}.
Accor\-ding to \cite[Prop.~2.4, Thm.~3.4, Thm.~3.7]{Ay2}, they
commute with the functors $f^*$, $f_*$, $f_!$, $f^!$, provided the latter
are applied to constructible objects (note that commutation holds without
this restriction for the two functors $f^*$ and $f_!$). 
In particular, they commute with Tate twists.
In \cite[Ex.~17.1.7]{CD},
it is shown how to obtain from the $R_{\eta,Z}$ a family of exact functors
with analogous properties, and which we denote by the same symbols
\[
R_{\eta,Z} : \DBcZM \longto D^b_c(Z) \; ,
\]     
where the right hand side denotes the full triangulated sub-category
of $D(Z)$ of classes of bounded complexes with constructible
cohomology objects. The construction can be imitated to obtain  
$F$-linear versions of the Betti realization. Composing with the 
base change \emph{via} $\Spec k \into \BB$, we finally obtain a family
of exact tensor functors
\[
R_{\eta, X} : \DBcFXM \longto D^b_c(X_k)_F \; ,
\]  
still referred to as the Betti realization,   
and indexed by schemes $X$, whose base change
$X_k := X \times_\BB \Spec k$ is quasi-projective over $k$. The functors $R_{\eta, X}$
are symmetric monoidal; in particular, they
respect the unit objects. They commute with the functors 
$f^*$, $f_*$, $f_!$, $f^!$ since $\Spec k \into \BB$ is a projective limit
of open immersions (use \cite[Prop.~14.3.1]{CD}),
\item[(ii)] $k$ is of characteristic zero, and $\ell$ is a prime. 
The \emph{$\ell$-adic realization} is defined in 
\cite[Sect.~7.2, see in part.\ Rem.~7.2.25]{CD2}. 
It is a family of covariant exact functors
\[
R_{\ell,Z}: \DBcZM \longto D^b_c(Z) \; ,
\] 
indexed by $k$-schemes $Z$ of finite dimension. Its target is 
the bounded ``derived category'' of
constructible $\BQ_\ell$-sheaves on $Z$ \cite[Sect.~6]{E}.
The functors $R_{\ell,Z}$ 
are symmetric monoidal, and they
commute with the functors $f^*$, $f_*$, $f_!$, $f^!$ \cite[Thm.~7.2.24]{CD2}. 
In parti\-cular, they commute with Tate twists.
The construction can be imitated to obtain  
$F$-linear versions of the $\ell$-adic realization. Composing with the 
base change \emph{via} $\Spec k \into \BB$, we finally obtain a family
of exact tensor functors
\[
R_{\ell, X} : \DBcFXM \longto D^b_c(X_k)_F \; ,
\]  
still referred to as the $\ell$-adic realization. The $R_{\ell, X}$
are symmetric monoidal, and they
commute with the functors $f^*$, $f_*$, $f_!$, $f^!$.
\end{enumerate}

In both settings, the categories $D^b_c(X_k)_F$ 
are equipped with a perverse $t$-structure;
write $D^{t=0}(X_k)$ for its heart,  
$H^n : D^b_c(X_k)_F \to D^{t=0}(X_k)$, $n \in \BZ$, 
for the cohomology functors, and
\[
H^* \! R_X : = \bigr( H^n \! R_X \bigl)_{n \in \BZ} : 
\DBcFXM \longto \Gr_{\BZ} D^{t=0}(X_k) 
\]
for the collection of all cohomology functors, preceded by the
realization $R_X$.
Here, we denote by $\Gr_{\BZ} D^{t=0}(X_k)$ 
the $\BZ$-graded category associated
to the heart $D^{t=0}(X_k)$. We shall often refer to $H^* \! R_X$ as 
the \emph{cohomological realization}.

\begin{Rem}
If $k$ is finitely generated over $\BQ$, then
the $\ell$-adic realization is a realization ``with weights'' 
as the action of local Frobenii allows for a notion of purity. 
By contrast, there are no intrinsic weights on the target
category of the Betti realization. The ideal solution would be to replace  it
by the \emph{Hodge theoretical realization} 
\[
R_{{\bf H} , \eta, X} : \DBcFXM \longto 
D^b \bigl( \MHM_\BQ  (X \times_\eta \Spec \BC) \otimes_\BQ F \bigr) 
\] 
to the bounded derived category 
of \emph{algebraic mixed Hodge modules} on $X \times_\eta \Spec \BC$ 
\cite[Sect.~4.2]{Sa}.
\end{Rem}

Let 
\[
R_{\bullet} : \DBcFbM \longto D^b_c(\bullet_k)_F 
\]  
be one of the two realizations considered above (Betti or $\ell$-adic).

\begin{Prop} \label{1Lpre}
Let $Y = Y(\Phi) = \coprod_{\varphi \in \Phi} \Yp$ be a good stratification
with nilregular strata,
such that the generic points of all $\Yp$ lie over $\Spec k \into \BB$. 
In the context of the Betti realization, assume that $Y_k$
is quasi-projective over $\Spec k$. Then the restriction of the
cohomological realization functor on $Y$ to $\DBcFYPAbM$,
\[
H^* \! R_{Y} : \DBcFYPAbM \longto \Gr_{\BZ} D^{t=0}(Y_k)
\]
satisfies assumptions~(1), (2), (4) and (5)
of Theorem~\ref{0G}. If $R_{\bullet}$ is the $\ell$-adic realization,
then 
\[
H^* \! R_{\ell,Y} : \DBcFYPAbM \longto \Gr_{\BZ} D^{t=0}(Y_k)
\]
also satisfies assumption~(3) of Theorem~\ref{0G}.
\end{Prop}

\begin{Proof}
Boundedness of the weight structure on $\DBcFYPAbM$ is
Theorem~\ref{1G}~(b). By Theorem~\ref{1I}, its heart 
$\CHFYPAbM$ is semi-primary and pseudo-Abelian. 
According to \cite[Cor.~7.13]{W9} (see Lemma~\ref{1lem}),
the restriction of $H^* \! R_{Y}$ to 
$\DBcFYPAbM$ maps the radical to zero. 

Thus, assumptions~(1), (2) and (4)
of Theorem~\ref{0G} are met.
Let us check assumption~(5), \emph{i.e.}, let us show that the zero
motive is the only Chow motive in $\CHFYPAbM$
whose realization is zero. 

First, given Theorem~\ref{1K}, we may assume that $\BB$ equals the generic point
$\Spec k$. Thus, we have $Y_k = Y$.

Second, by Definition~\ref{1F}~(b), and by \cite[Cor.~4.10~(b)]{W9},
the triangulated category $\DBcFYPAbM$
is obtained by successive gluing over the strata
$\Yp$ of triangulated sub-categories $\DBcFYppAbM$ of $\DBcFYpAbM$. 
The $\DBcFYppAbM$ inherit the weight structure (Theorem~\ref{1G}) from $\DBcFYpAbM$,
and according to Theorem~\ref{1I}, their hearts are semi-primary and pseudo-Abelian.
Thus, the abstract theory of \emph{intermediate extensions} can be applied: according to
\cite[Summ.~2.12]{W9}, any object of $\CHFYPAbM$ is a direct sum of Chow motives of the
form $j_{\varphi,!*} M_\varphi$, for certain Chow motives 
$M \in \CHFYpAbM$. Here, $j_{\varphi,!*}$ denotes the
intermediate extension \cite[Def.~2.10]{W9} associated to the immersion 
$j_{\varphi} : \Yp \into Y$. If the cohomological realization 
$H^* \! R_{\Yp}(M) \in \Gr_{\BZ} D^{t=0}(\Yp)$ is zero, then so are the
cohomological realizations of all $j_{\varphi,!*} M_\varphi$, hence of all $M_\varphi$
since $R_{\bullet}$ is compatible with inverse images. 
Thus, we may assume that the 
stratification $\Phi$ consists of a single (nilregular)
stratum: $Y(\Phi) = \Yp$.

According to Lemma~\ref{1lem}, \cite[Prop.~5.5]{W9} and Proposition~\ref{1Ipre},
any Chow motive $M$ in $\CHFYpAbM$ is smooth and finite dimensional.
The same is therefore true for its pull-back $M_\xi$
to any generic point $\xi$ of $\Yp$. 

So if we assume the
cohomological realization 
$H^* \! R_{\Yp}(M) \in \Gr_{\BZ} D^{t=0}(\Yp)$ to be zero,
then the realization of any $M_\xi$ is zero, again
since $R_{\bullet}$ is compatible with inverse images.
Therefore \cite[Cor.~7.3]{Ki}, all $M_\xi$ are zero. 
But according to \cite[pp.~54--55]{O'S}, this implies that $M$ is zero.

If $R_{Y}$ is the $\ell$-adic
realization, then hypothesis~(3) of Theorem~\ref{0G}
is also met, \emph{i.e.}, the morphisms in 
the image of $H^* \! R_{\ell,Y}$ are strict
with respect to the weight filtration
\cite[Thm.~2.5.4~(II)~(1), Prop.~1.3.2~(II)~(2)]{Bo2}.
Note that $Y_k$ 
is of finite type over $k$, hence \emph{very reasonsable}
in the sense of \cite[Def.~2.1.1~(4)]{Bo2}.
\end{Proof}

Here is our second result on conservativity in the motivic context.

\begin{Thm} \label{1L}
(a)~Let $Y = Y(\Phi) = \coprod_{\varphi \in \Phi} \Yp$ be a good stratification
with nilregular strata,
such that the generic points of all $\Yp$ lie over $\Spec k \into \BB$. 
In the context of the Betti realization, assume that $Y_k$
is quasi-projective over $\Spec k$. Then the restriction of the
realization functor on $Y$ to $\DBcFYPAbM$,
\[
R_{Y} : \DBcFYPAbM \longto D^b_c(Y_k)_F \; ,
\]
is conservative. \\[0.1cm]
(b)~Assume that $\BB = \Spec k$. Let $Y = Y(\Phi) = \coprod_{\varphi \in \Phi} \Yp$ 
be a good stratification with nilregular strata. 
In the context of the Betti realization, assume that $Y$
is quasi-projective over $\Spec k$. Then the restriction of the
realization functor on $Y$ to $\DBcFYPAbM$,
\[
R_{Y} : \DBcFYPAbM \longto D^b_c(Y)_F \; ,
\]
is conservative.
\end{Thm}

\begin{Proof}
(b) is a special case of (a) (put $\BB = \Spec k$).

For the $\ell$-adic realization, Theorem~\ref{0G} can be applied
directly, thanks to Proposition~\ref{1Lpre}.

Let us treat the case when $R_{\bullet}$ is the Betti realization.
Given our present state of knowledge, we may suppose, but do not know
the analogue of hypo\-the\-sis~(3) of Theorem~\ref{0G} to hold; 
therefore, we need an alternative approach.

As earlier (Theorem~\ref{1K} and its proof), 
we may assume that $\BB$ equals the generic point
$\Spec k$, and that the 
stratification of $Y(\Phi)$ consists of a single (nilregular)
stratum: $Y(\Phi) = \Yp$.

Now recall from the proof of Theorem~\ref{0G} that hypothesis~(3)
is only used \emph{via} Lemma~\ref{0F}.
The idea therefore consists in deducing Lemma~\ref{0F} from 
the little we know. By
Lemma~\ref{1lem}, \cite[Prop.~5.5]{W9} and Proposition~\ref{1Ipre},
the category $\CHFYpAbM$
consists of Chow motives which are smooth over $\Yp$. As it generates
the triangulated category $\DBcFYpAbM$ 
(Theorem~\ref{1G}~(b)), the cohomological realization of any object
of the latter gives perverse sheaves which are actually local systems
(up to a shift). 
This holds in particular for the objects $X$ and $Z$
occurring in Lemma~\ref{0F}; therefore, the effect of the cohomological realization of
a morphism between them can be read off the restriction of the latter to
the generic points of $\Yp$, where comparison with the $\ell$-adic realization
is available.
\end{Proof}

The third and main result on conservativity reads as follows;
it generalizes the $\ell$-adic version of \cite[Thm.~1.13]{W8}.

\begin{Thm} \label{1M}
Assume $k$ to be of characteristic zero, and let $\ell$ a prime.
Let $Y = Y(\Phi) = \coprod_{\varphi \in \Phi} \Yp$ be a good stratification
with nilregular strata,
such that the generic points of all $\Yp$ lie over $\Spec k \into \BB$. 
Then the $\ell$-adic realization $R_{\ell,Y}$ respects and 
detects the weight structure on $\DBcFYPAbM$. More precisely,
let $M \in \DBcFYPAbM$, and $\alpha \le \beta$ two integers. \\[0.1cm]
(a)~$M$ lies in the heart $\CHFYPAbM$ of $w$ 
if and only if 
the $n$-th perverse cohomology object 
$H^n \! R_{\ell,Y}(M) \in D^{t=0}(Y_k)$ of $R_{\ell,Y}(M)$ 
is pure of weight $n$, 
for all $n \in \BZ$. \\[0.1cm]
(b)~$M$ lies in $DM_{\text{\cyrb},c,\Phi}^{Ab}(Y)_{F,w \le \alpha}$ if and only if 
$H^n \! R_{\ell,Y}(M)$ is of weights $\le n + \alpha$, 
for all $n \in \BZ$. \\[0.1cm] 
(c)~$M$ lies in $DM_{\text{\cyrb},c,\Phi}^{Ab}(Y)_{F,w \ge \beta}$ if and only if 
$H^n \! R_{\ell,Y}(M)$ is of weights $\ge n + \beta$, 
for all $n \in \BZ$. \\[0.1cm]
(d)~$M$ is without weights $\alpha,\alpha+1,\ldots,\beta$ if and only if 
$H^n \! R_{\ell,Y}(M)$ is without weights 
$n + \alpha,n + \alpha+1,\ldots,n + \beta$, for all $n \in \BZ$. 
\end{Thm}

\begin{Proof}
According to Proposition~\ref{1Lpre},
the assumptions of Theorem~\ref{0H} are satisfied for 
$\CH = H^0 \! R_{\ell,Y}$. 
\end{Proof}

\begin{Rem}
If $k$ is finitely generated over $\BQ$, then there is an intrinsic notion
of weights on those objects of the heart $D^{t=0}(Y_k)$ 
of the perverse $t$-structure on $D^b_c(Y_k)_F$, which are in the
image of the cohomological realization \cite[Prop.~2.5.1~(II)]{Bo2}.

In general, the weights of $H^* \! R_{\ell,Y}(M)$ are by definition
those induced by the weight filtration of the functor
$H^* \! R_{\ell,Y}$ as considered in the previous section 
(these coincide with the
above when $k$ is finitely generated over $\BQ$).
\end{Rem}

\begin{Rem}
The analogue of Theorem~\ref{1M} should hold for the Betti, and the Hodge
theoretical realization. In the absence of the latter, 
and/or a general comparison statement between the Betti and the $\ell$-adic
realization, the problem for the Betti realization
consists in the verification of assumption~(3) of Theorem~\ref{0G}
(which \emph{directly} enters the proof of Theorem~\ref{0H}).
Contrary to the proof of Theorem~\ref{1L},
reduction to a statement on individual strata of $Y(\Phi)$ does not seem
to work. 
\end{Rem}

Let us spell out the special case $\BB = \Spec k$ of Theorem~\ref{1M}.

\begin{Cor} 
Let $k$ be a field of characteristic zero, $\ell$ a prime,
$Y$ a scheme over $\BB = \Spec k$, and
$Y = Y(\Phi) = \coprod_{\varphi \in \Phi} \Yp$ a good stratification
with nilregular strata. 
Let $M \in \DBcFYPAbM$, and $\alpha \le \beta$ two integers. \\[0.1cm]
(a)~$M$ lies in the heart $\CHFYPAbM$ of $w$ 
if and only if 
the $n$-th perverse cohomology object 
$H^n \! R_{\ell,Y}(M) \in D^{t=0}(Y)$ of $R_{\ell,Y}(M)$ 
is pure of weight $n$, 
for all $n \in \BZ$. \\[0.1cm]
(b)~$M$ lies in $DM_{\text{\cyrb},c,\Phi}^{Ab}(Y)_{F,w \le \alpha}$ if and only if 
$H^n \! R_{\ell,Y}(M)$ is of weights $\le n + \alpha$, 
for all $n \in \BZ$. \\[0.1cm] 
(c)~$M$ lies in $DM_{\text{\cyrb},c,\Phi}^{Ab}(Y)_{F,w \ge \beta}$ if and only if 
$H^n \! R_{\ell,Y}(M)$ is of weights $\ge n + \beta$, 
for all $n \in \BZ$. \\[0.1cm]
(d)~$M$ is without weights $\alpha,\alpha+1,\ldots,\beta$ if and only if 
$H^n \! R_{\ell,Y}(M)$ is without weights 
$n + \alpha,n + \alpha+1,\ldots,n + \beta$, for all $n \in \BZ$. 
\end{Cor}


\bigskip

%
%

\end{document}